# LIKELIHOOD APPROACH FOR MARGINAL PROPORTIONAL HAZARDS REGRESSION IN THE PRESENCE OF DEPENDENT CENSORING[1]

By Donglin Zeng

*University of North Carolina at Chapel Hill*

In many public health problems, an important goal is to identify the effect of some treatment/intervention on the risk of failure for the whole population. A marginal proportional hazards regression model is often used to analyze such an effect. When dependent censoring is explained by many auxiliary covariates, we utilize two working models to condense high-dimensional covariates to achieve dimension reduction. Then the estimator of the treatment effect is obtained by maximizing a pseudo-likelihood function over a sieve space. Such an estimator is shown to be consistent and asymptotically normal when either of the two working models is correct; additionally, when both working models are correct, its asymptotic variance is the same as the semiparametric efficiency bound.

**1. Introduction.** In many public health problems, an important goal is to study the effect of some treatment or intervention on the risk of failure. A commonly used model to analyze such an effect is via the proportional hazards regression model:

$$(1.1) \qquad h_{T|V}(t|v) = \lambda(t) e^{\alpha v},$$

where $V$ denotes the measurement of treatment, $T$ denotes failure time and $h_{T|V}(t|v)$ denotes the hazard rate function of $T$ given $V$. In the model (1.1), $\lambda(t)$ is an unknown baseline hazard rate function and $\alpha$ is an unknown parameter describing the effect of $V$. A marginal regression model such as (1.1) is often useful in public health problems, since in that field the scientific goal is to identify the effect of treatment for the whole population regardless of heterogeneity within the population; in other words, we would not adjust

Received January 2002; revised May 2004.
[1]Supported in part by NIDA Grant 2P50 DA 10075 and NSF Grant DMS-98-02885.
*AMS 2000 subject classifications.* Primary 62G07; secondary 62F12.
*Key words and phrases.* Semiparametric inference, dimension reduction, B-spline, double robustness.

for other covariates in the regression model (1.1) even if such covariates are measured at the same time of data collection. Some other reasons why additional covariates would not be adjusted for in the regression model for epidemiologic studies can be seen in Robins, Rotnitzky and Zhao (1994).

Dependent right-censoring is common in failure time data, where subjects may drop out or be censored during the studies. The censorship can be caused by many factors, such as the feeling of patients about participation in the studies, the social supports for patients, patients' accessibility to the studies, biological information of patients, and so on. In practice, when a large amount of such information is collected, it is safe to assume that the dependence between the failure time and the censoring time is fully explained by all the collected covariates. In mathematical notation, if we denote $C$ as censoring time and denote $X$ as other auxiliary covariates besides $V$, then we assume that, conditional on $X$ and $V$, $T$ and $C$ are independent.

Suppose $n$ i.i.d. right-censored observations are available and we denote them as $(Y_i = T_i \wedge C_i, R_i = I(T_i \leq C_i), X_i, V_i)$, $i = 1, \ldots, n$. Our goal is to estimate the treatment effect $\alpha$ in the model (1.1). It is well known that, in the presence of dependent censoring, simply performing the Cox regression using $V$ as covariates gives an inconsistent estimate. In order to adjust for dependent censoring, one intuitive approach tends to estimate the distribution of $T$ given $(X, V)$ either nonparametrically or semiparametrically. However, two weak points can restrain the use of this approach: one is that nonparametric estimation is not feasible with moderate samples when $X$ has more than three dimensions, which is known as the curse of the dimensionality; the other is that many semiparametric models of $T$ given $(X, V)$ are generally not compatible with (1.1) while the latter, as indicated above, is of main scientific interest. Recently, an estimating equation approach was proposed by Robins, Rotnitzky and Zhao (1994) and was successfully applied to missing longitudinal data; however, to our knowledge, such an approach has not been applied to regression problems for survival data, except for a brief discussion in Robins, Rotnitzky and Zhao (1994). Furthermore, the implementation of the estimating equation approach relies on the derivation of the efficient score function for $\alpha$, which is implicit and difficult for the model (1.1).

In this paper, we propose a likelihood-based approach to estimate the parameters in the marginal proportional hazards model (1.1). The ideas of handling dependent censoring are similar to those in one of our previous papers [Zeng (2004)]. Briefly, we first condense the high-dimensional covariates $(X, V)$ by utilizing two working models for the distribution of $T$ given $(X, V)$ and the distribution of $C$ given $(X, V)$. Then an estimate for the coefficient $\alpha$ in (1.1) is obtained by maximizing a pseudo-likelihood function of a reduced datum, which consists of the observed event times, the censoring status, the treatment and the condensed covariates. In the maximization, the nuisance parameters for $\alpha$ are profiled out over a sieve space consisting of B-splines. At the end of this paper we demonstrate that the estimator for $\alpha$ has the following properties: if either of the two working models is correct, the estimator is consistent and asymptotically normal; if both working models are correct, the estimator's asymptotic variance attains the semiparametric efficiency bound. The first property is called double robustness by Robins, Rotnitzky and van der Laan (2000). The details of the proofs are given in the Appendix.

**2. Estimation.** For convenience, we denote $f_{Z_1|Z_2}(\cdot|\cdot)$ as the conditional density of a random vector $Z_1$ given another random vector $Z_2$ and abbreviate $(X^T, V)^T$ as $W$.

2.1. *Estimation procedure.* First we utilize two working models for the distribution of $T$ and $C$ given $W$.

WORKING MODEL 1. We tentatively assume that $T$ is independent of $X$ given $V$ so $f_{T|W}(y|w) = f_{T|V}(y|v)$.

WORKING MODEL 2. We tentatively assume that the model of $C$ given $W$ is a proportional hazards model, that is, $h_{C|W}(y|w) = h_c(y) e^{\gamma^T W}$ for an unknown vector $\gamma$ and an unknown baseline hazard rate function $h_c(\cdot)$.

REMARK. In fact, any model can be used for Working Model 1 and there are two reasons for us to choose the current form: first, this is a simple one to work with; second, our later results show that, to ensure our proposed estimator is more likely to be consistent, such a working model has to satisfy the constrained form in (1.1). Obviously, the current Working Model 1 is the most convenient choice.

To illustrate the estimation procedure, we suppose that either working model is correct and that $\gamma$ is a known constant. We let $[a(\gamma, v), b(\gamma, v)]$ be the support of the conditional distribution of $\gamma^T W$ given $V = v$ and define

$$U(\gamma) = [\gamma^T W - a(\gamma, V)]/[b(\gamma, V) - a(\gamma, V)]$$

for fixed $\gamma$. Then the conditional distribution of $U(\gamma)$ given $V$ has support $[0, 1]$. As shown in Lemma 3.1 of Zeng (2004), $T$ and $C$ are independent given $(U(\gamma), V)$ when either working model is correct; in other words, the dependence between $T$ and $C$ can be fully explained by the two-dimensional condensed information $(U(\gamma), V)$. We replace the observed statistics $W$ with $(U(\gamma), V)$ and obtain reduced data $(Y_i, R_i, U_i, V_i)$, where $U_i = \gamma^T W_i$, $i =$ $1, \ldots, n$. Therefore, the observed likelihood function of the reduced data concerning the joint distribution of $(T, U(\gamma))$ given $V$ is

$$\prod_{i=1}^{n} \left\{ [f_{T|U(\gamma),V}(Y_i|U_i,V_i)]^{R_i} \left[ \int_{Y_i}^{\infty} f_{T|U(\gamma),V}(s|U_i,V_i)\, ds \right]^{1-R_i} f_{U(\gamma)|V}(U_i|V_i) \right\}.$$

In order to absorb the marginal model (1.1) into the observed likelihood function, a natural reparameterization is to use the conditional density of $U(\gamma)$ given $T$ and $V$ and the conditional density of $T$ given $V$ as the new parameters. The latter contains the parameters $\lambda(y)$ and $\alpha$. However, since $T$ is only observable in $[0, \tau)$ where $\tau$ is the end time of the study, the conditional density of $U(\gamma)$ given $T$ and $V$ is not identifiable for $T \geq \tau$. Therefore, we introduce a modified new variable $\widetilde{T} = TI(T < \tau) + \tau I(T \geq \tau)$; that is, $\widetilde{T}$ is the same as $T$ if $T$ is observed within the study time frame and $\widetilde{T}$ is equal to $\tau$ if $T$ is out of the observable range. Then it is easy to calculate the density function for $\widetilde{T}$ given $V = v$ as

$$I(t < \tau)\lambda(t)e^{\alpha v}\exp\{-e^{\alpha v}\Lambda(t)\} + \delta(t = \tau)\exp\{-e^{\alpha v}\Lambda(\tau)\},$$

where $\delta(\cdot)$ is the Dirac function. Moreover, we denote $f_{U(\gamma)}(\cdot|y,v)$ as the conditional density of $U(\gamma)$ given $\widetilde{T} = y$ and $V = v$ for $y \in (0, \tau)$ and denote $g_{U(\gamma)}(u|\tau, v)$ as the conditional density of $U(\gamma)$ given $\widetilde{T} = \tau$ and $V = v$. Thus, $f_{U(\gamma)}(\cdot|y,v)$ is the same as the conditional density of $U(\gamma)$ given $T = y$ and $V = v$ and $g_{U(\gamma)}(\cdot|\tau, v)$ is the same as the conditional density of $U(\gamma)$ given $T \geq \tau$ and $V = v$. Since the observed data are equivalent to $(U_i, V_i, R_i = I(\widetilde{T}_i \leq C_i), Y_i = \widetilde{T}_i \wedge C_i)$, in terms of the new parameters $(\alpha, \lambda(y), f_{U(\gamma)}(u|y,v), g_{U(\gamma)}(u|\tau, v))$ the observed likelihood function can be written as

$$\prod_{i=1}^{n} \Big\{ [\exp\{-e^{\alpha V_i}\Lambda(Y_i)\}e^{\alpha V_i}\lambda(Y_i)f_{U(\gamma)}(U_i|Y_i,V_i)]^{R_i}$$

(2.1)
$$\times \left[ \int_{Y_i}^{\tau-} \exp\{-e^{\alpha V_i}\Lambda(s)\}e^{\alpha V_i}\lambda(s)f_{U(\gamma)}(U_i|s,V_i)\, ds \right.$$

$$\left. + \exp\{-e^{\alpha V_i}\Lambda(\tau)\}g_{U(\gamma)}(U_i|\tau,V_i) \right]^{1-R_i} \Big\}.$$

Clearly, all the parameters are distinct and identifiable.

Finally, we maximize the function (2.1) over a sieve space of the parameters $(\alpha, \lambda(y), f_{U(\gamma)}(u|y,v), g_{U(\gamma)}(u|\tau, v))$ for some estimate of $\gamma$. In the following sections, we describe how to obtain an estimate of $\gamma$ and how to construct a sieve space for the parameters.

2.2. *An estimate for $\gamma$.* We estimate $\gamma$ by performing the proportional hazards regression using the censored observations. That is, we maximize the following pseudo-partial likelihood function for $\gamma$:

$$\prod_{i=1}^{n}\left[\frac{e^{\gamma^T W_i}}{\sum_{Y_j \geq Y_i} e^{\gamma^T W_j}}\right]^{1-R_i}.$$

The estimator for $\gamma$ is denoted as $\hat{\gamma}_n$. As shown in Theorem 3.1 of Zeng (2004), under some regular conditions $\hat{\gamma}_n$ should converge to a constant $\gamma^*$ almost surely and $\sqrt{n}(\hat{\gamma}_n - \gamma^*)$ has an asymptotically linear expansion with its influence function denoted by $S(Y, R, W; \gamma^*)$.

2.3. *Sieve space $\mathcal{S}_n$ for the parameters $(\alpha, \lambda(y), f_{U(\gamma)}(u|y,v), g_{U(\gamma)}(u|\tau,v))$.* We propose a sieve space consisting of B-splines for $f_{U(\gamma)}(u|y,v)$, $g_{U(\gamma)}(u|\tau,v)$ and $\lambda(y)$ in maximizing (2.1). We suppose that $0 \leq V \leq 1$ and that $|\alpha| \leq M$ for a known constant $M$.

We reparameterize $(f_{U(\gamma)}(u|y,v), g_{U(\gamma)}(u|\tau,v), \lambda(y))$ by introducing

$$f_{U(\gamma)}(u|y,v) = \frac{\exp\{\eta_1(u,y,v)\}}{\int_0^1 \exp\{\eta_1(u,y,v)\}\,du},$$

$$g_{U(\gamma)}(u|\tau,v) = \frac{\exp\{\eta_2(u,v)\}}{\int_0^1 \exp\{\eta_2(u,v)\}\,du},$$

and $\lambda(y) = \exp\{\xi(y)\}$, where $\eta_1(u,y,v)$ and $\eta_2(u,v)$ satisfy that $\eta_1(0,y,v) = 0$, $\eta_2(0,v) = 0$. After the reparameterization, the new parameters are $(\alpha, \xi(y), \eta_1(u,y,v), \eta_2(u,v))$ in which $0 \leq u,v \leq 1$, $0 \leq y \leq \tau$. A sieve space consisting of B-splines is defined for these new parameters as follows: First, we obtain an extended partition with equal length $1/K_n$ for the interval $[0,1]$:

$$\Delta_e = \{s_{-m} = \cdots = s_{-1} = 0 = s_0 < s_1 < \cdots < s_{K_n} = 1 = \cdots = s_{K_n+m}\},$$

where $m$ (independent of $n$) and $K_n$ are two integers to be chosen later. Let $\{N_j^m(s)\}_{j=1}^{K_n+m}$ be a normalized B-spline basis associated with $\Delta_e$ [cf. Schumaker (1981)]. Then the sieve space for the parameters $(\alpha, \xi(y), \eta_1(u,y,v), \eta_2(u,v))$ is defined as

(2.2)
$$\begin{aligned}&\mathcal{S}_n(m, K_n, M_n)\\&= \Bigg\{(\alpha, \xi(y), \eta_1(u,y,v), \eta_2(u,v)) : |\alpha| \leq M,\\&\quad \eta_1(u,y,v) = \sum_{i_1,i_2,i_3=1}^{m+K_n} \eta^1_{i_1,i_2,i_3} N_{i_1}^m(u) N_{i_2}^m(y/\tau) N_{i_3}^m(v),\end{aligned}$$

$$\eta_2(u,v) = \sum_{i_1,i_2=1}^{m+K_n} \eta_{i_1,i_2}^2 N_{i_1}^m(u) N_{i_2}^m(v), \ \xi(y) = \sum_{i=1}^{m+K_n} \xi_i N_i^m(y/\tau),$$

$$\sum_{i_1,i_2,i_3=1}^{m+K_n} |\eta_{i_1,i_2,i_3}^1| \le M_n, \ \sum_{i_1,i_2=1}^{m+K_n} |\eta_{i_1,i_2}^2| \le M_n, \ \sum_{i=1}^{m+K_n} |\xi_i| \le M_n,$$

$$\left. \sum_{i_1=1}^{m+K_n} \eta_{i_1,i_2,i_3}^1 N_{i_1}^m(0) = 0, \ \sum_{i_1=1}^{m+K_n} \eta_{i_1,i_2}^2 N_{i_1}^m(0) = 0 \right\}.$$

In other words, we use a finite linear combination of the B-splines to approximate each nonparametric function. The use of the last two constraints in the conditions of the sieve space ensures that $\eta_1(0,y,v) = 0$ and $\eta_2(0,v) = 0$. The constants $M_n$ and $K_n$ depend on $n$ and will be chosen later.

2.4. *Maximization.* Let $\mathbf{P}_n, \mathbf{P}$ denote the empirical measure and the true probability measure of $(Y, R, W)$, respectively, and let $\hat{U} = [\hat{\gamma}_n^T W - a(\hat{\gamma}_n, V)]/[b(\hat{\gamma}_n, V) - a(\hat{\gamma}_n, V)]$. We maximize the function

$$\mathbf{P}_n \left\{ R \log \left[ \exp \left\{ -\int_0^Y e^{\xi(s)+\alpha V} \, ds \right\} e^{\xi(Y)+\alpha V} \frac{\exp\{\eta_1(\hat{U}, Y, V)\}}{\int_0^1 \exp\{\eta_1(u, Y, V)\} \, du} \right] \right\}$$

$$+ \mathbf{P}_n \left\{ (1-R) \log \left[ \int_Y^{\tau-} \exp \left\{ -\int_0^s e^{\xi(s')+\alpha V} \, ds' \right\} e^{\xi(s)+\alpha V} \right. \right.$$

(2.3)
$$\times \frac{\exp\{\eta_1(\hat{U}, s, V)\}}{\int_0^1 \exp\{\eta_1(u, s, V)\} \, du} \, ds$$

$$\left. \left. + \exp \left\{ -\int_0^\tau e^{\xi(s)+\alpha V} \, ds \right\} \frac{\exp\{\eta_2(\hat{U}, V)\}}{\int_0^1 \exp\{\eta_2(u, V)\} \, du} \right] \right\}$$

over the sieve space $\mathcal{S}_n(m, K_n, M_n)$. One possible choice of $(m, K_n, M_n)$ is $(k+2, \tilde{M}n^\beta, \tilde{M}\sqrt{\log n})$ for some given constant $\tilde{M}$, a known integer $k \ge 11$ and a constant $\beta$ satisfying $\frac{1}{2k} < \beta < \frac{3}{4k+9}$.

It will be shown later that $|a(\hat{\gamma}_n, v) - b(\hat{\gamma}_n, v)|$ has a positive limit for any $v$ with probability 1. Then the arguments of the maximum exist since we are maximizing the function over a compact set in a finite-dimensional space. However, the solution itself may not be unique. We simply select any one of these maximizers and denote it as $(\hat{\alpha}_n, \hat{\xi}_n(y), \hat{\eta}_{1n}(u,y,v), \hat{\eta}_{2n}(u,v))$. Respectively, we obtain the estimators $\hat{\alpha}_n = \hat{\alpha}_n$, $\hat{\lambda}_n(y) = \exp\{\hat{\xi}_n(y)\}$ and

$$\hat{f}_n(u|y,v) = \frac{\exp\{\hat{\eta}_{1n}(u,y,v)\}}{\int_0^1 \exp\{\hat{\eta}_{1n}(u,y,v)\} \, du}, \qquad \hat{g}_n(u|\tau,v) = \frac{\exp\{\hat{\eta}_{2n}(u,v)\}}{\int_0^1 \exp\{\hat{\eta}_{2n}(u,v)\} \, du}.$$

Computationally, many constrained optimization algorithms such as the quasi-Newton method, combined with the use of either a penalty or a barrier function, can be applied to find the arguments of the maximization.

**3. Asymptotic results.** We provide the main results in this section. Especially, the consistency and asymptotic distribution for $\hat{\alpha}_n$ are derived. The proofs for all the theorems are given in Section 6.

3.1. *Assumptions.* In addition to the assumption that $T$ and $C$ are independent given $W$, we need the following conditions.

ASSUMPTION A1. $V$ has support in $[0,1]$ and $X$ has bounded support in $R^d$ where $d$ is the dimension of $X$. Moreover, if there exist a constant $c_0$ and a constant vector $\tilde{\gamma}$ such that $\tilde{\gamma}^T W = c_0$ almost surely, then $c_0 = 0$ and $\tilde{\gamma} = 0$.

ASSUMPTION A2. With probability 1, there exists a positive constant $\theta_0$ such that $P(C \geq \tau|W) = P(C = \tau|W) \geq \theta_0$ and $P(T > \tau|W) \geq \theta_0$. That is, at least some subjects do not fail at the end time $\tau$ and by definition they are considered to be right-censored at $\tau$.

ASSUMPTION A3. For a known integer $k \geq 11$, the conditional density of $X$ given $\widetilde{T}$ and $V$, denoted as $f_{X|\widetilde{T},V}$, and the true baseline hazard rate, $\lambda_0(y)$, satisfy

$$\log f_{X|\widetilde{T},V}(x|y,v) \in W^{k+4,2}(R^{d+2}), \qquad \log \lambda_0(y) \in W^{k+4,2}(R),$$

after appropriate extension to the whole space. Here, $W^{k+4,2}(R^l)$ is a Sobolev space consisting of the functions with $(k+4)$th derivatives in $L_2(R^l)$. Furthermore, we assume that

$$\log f_{C|W}(y|w) \in W^{k+4,2}((0,\tau) \times R^{d+1}),$$
$$\log P(C = \tau|W = w) \in W^{k+4,2}(R^{d+1}).$$

ASSUMPTION A4. There exists a known constant $M$ such that the true treatment effect $\alpha_0$ satisfies $|\alpha_0| \leq M$. Moreover, the equation

$$\mathbf{P}[(1-R)W] = \mathbf{P}\left\{(1-R)\left[\frac{\mathbf{P}[I_{Y \geq y} W e^{\gamma^T W}]}{\mathbf{P}[I_{Y \geq y} e^{\gamma^T W}]}\right]\bigg|_{y=Y}\right\}$$

has a unique solution $\gamma^*$ in $[-M, M]^{d+1}$. In addition, for any $\gamma$ in a small neighborhood $\mathcal{O}$ of $\gamma^*$, the conditional distribution of $\gamma^T W$ given $V = v$ has support $[a(\gamma, v), b(\gamma, v)]$ satisfying: both the function $a(\cdot)$ and the function $b(\cdot)$ are two known functions and they are continuously differentiable with respect to $\gamma$; as functions of $v$, $a(\gamma, v)$ and $b(\gamma, v)$ belong to $W^{k+4,2}(R)$; $\min_{v \in [0,1], \gamma \in \mathcal{O}} |b(\gamma, v) - a(\gamma, v)| > 0$.

ASSUMPTION A5. $(M_n, K_n)$ satisfy $M_n, K_n \to \infty$ and

$$\frac{e^{13M_n}}{K_n^k} + \frac{e^{16M_n} K_n^{4k/3+3} (\log K_n)^2}{n} \to 0.$$

ASSUMPTION A6. $K_n$ satisfies $\sqrt{n} = o(K_n^{2k})$.

REMARK. Theorem 3.1 of Zeng (2004) showed that the asymptotic limit of $\hat{\gamma}_n$ is equal to $\gamma^*$ given by Assumption A4. It is also implied by Assumption A4 that one of the first $d$ components of $\gamma^*$ is nonzero. Thus, if we suppose the first component of $\gamma^* = (\gamma_1^*, \ldots, \gamma_d^*, \gamma_{d+1}^*)$ is not zero, then in terms of $f_{X|\widetilde{T},V}(x|y,v)$, the conditional density of $U^* = \gamma^{*T} W$ given $(\widetilde{T} = y, V = v)$, which is denoted as $f_{U^*}(u|y,v)$ for $\widetilde{T} < \tau$ and as $g_{U^*}(u|\tau,v)$ for $\widetilde{T} = \tau$, can be expressed as

$$\frac{\Delta(v)}{|\gamma_1^*|} \int f_{X|\widetilde{T},V} \bigg( \frac{u\Delta(v) + a(\gamma^*,v) - \sum_{i=2}^d \gamma_i^* x_i - \gamma_{d+1}^* v}{\gamma_1^*},$$
$$x_2, \ldots, x_d | y, v \bigg) dx_2 \cdots dx_d,$$

where $\Delta(v) = b(\gamma^*,v) - a(\gamma^*,v)$. Hence, Assumption A3 implies that $f_{U^*}(u|y,v)$ and $g_{U^*}(u|\tau,v)$ are bounded away from 0 and their $(k+4)$th derivatives are also $L_2$-integrable. Furthermore, by the embedding theorem in Sobolev space [cf. Adams (1975)], this gives that each of $\log f_{U^*|T,V}(u|y,v)$, $\log g_{U^*}(u|\tau,v)$, $\log f_{C|U^*,V}(y|u,v)$, $\log \lambda_0(y)$ is in $W^{k,\infty}$ space; that is, their $k$th derivatives are bounded essentially.

REMARK. Assumptions A5 and A6 determine the size of the sieve space in terms of the number of knots in the partition $(K_n)$ and the bounds of the sieve functions $(M_n)$. When $k \geq 11$, such $K_n$ satisfying both Assumptions A5 and A6 exists. For example, we can choose $K_n = n^\beta$, $\frac{1}{2k} < \beta < \frac{3}{4k+9}$. Additionally, the choice of $M_n$ can be of order $\sqrt{\log n}$.

Although all these assumptions guarantee the validity of the following arguments, they are not minimal assumptions.

3.2. *Asymptotic results.*

THEOREM 3.1 (Consistency of $\hat{\alpha}_n$). *Suppose that either of the two working models is true. Under Assumptions* A1–A5, *$\hat{\alpha}_n$ is a consistent estimator of the true coefficient $\alpha_0$.*

We can further obtain the consistency of the nuisance parameters in a Sobolev-norm.

THEOREM 3.2 (Consistency of nuisance parameters). *Suppose that either of the two working models is true. Under Assumptions* A1–A5,

$$\|\hat{\lambda}_n(Y) - \lambda_0(Y)\|_{W^{1,\infty}(P)} \xrightarrow{p} 0,$$

$$\|\hat{f}_n(U^*|Y,V) - f_{U^*}(U^*|Y,V)\|_{W^{1,\infty}(P)} \xrightarrow{p} 0,$$

$$\|\hat{g}_n(U^*|\tau,V) - g_{U^*}(U^*|\tau,V)\|_{W^{1,\infty}(P)} \xrightarrow{p} 0.$$

*Here* $\|h(U^*,Y,V)\|_{W^{1,\infty}(P)}$ *is defined as* $\|h(U^*,Y,V)\|_{L_\infty(P)} + \|\nabla h(U^*,Y,V)\|_{L_\infty(P)}$, *where* $P$ *is the probability measure given by* $(U^*,Y,V,R)$.

The result in Theorem 3.2 can help to obtain a useful convergence rate of the estimators in $L_2$-norm, which is stated in Theorem 3.3.

THEOREM 3.3 (Convergence rate). *Suppose that either of the two working models is true. Under Assumptions* A1–A5, *it holds that*

$$|\hat{\alpha}_n - \alpha_0|^2 + \|\hat{\lambda}_n(Y) - \lambda_0(Y)\|^2_{L_2(P)} \leq O_p\left(\frac{1}{K_n^{2k}}\right) + o_p\left(\frac{1}{\sqrt{n}}\right),$$

$$\|\hat{f}_n(U^*|Y,V) - f_{U^*}(U^*|Y,V)\|^2_{L_2(P)} \leq O_p\left(\frac{1}{K_n^{2k}}\right) + o_p\left(\frac{1}{\sqrt{n}}\right)$$

*and*

$$\|\hat{g}_n(U^*|\tau,V) - g_{U^*}(U^*|\tau,V)\|^2_{L_2(P)} \leq O_p\left(\frac{1}{K_n^{2k}}\right) + o_p\left(\frac{1}{\sqrt{n}}\right).$$

Finally, we derive the asymptotic distribution for $\sqrt{n}(\hat{\alpha}_n - \alpha_0)$.

THEOREM 3.4 (Asymptotic normality of $\hat{\alpha}_n$). *Under Assumptions* A1–A6, *when either of the two working models is correct,* $\sqrt{n}(\hat{\alpha}_n - \alpha_0)$ *is asymptotically normal. Furthermore, when both working models are correct, the asymptotic variance of* $\sqrt{n}(\hat{\alpha}_n - \alpha_0)$ *is the same as the semiparametric efficiency bound.*

3.3. *Variance estimation.* We propose the following steps to estimate the asymptotic variance of $\hat{\alpha}_n$ with no attempt to justify them rigorously. Our way is to directly estimate the influence function of $\hat{\alpha}_n$.

Define $O = (Y, R, W)$ and define $\psi$ as the nuisance parameters consisting of $(f_{U(\gamma)}(u|y,v), g_{U(\gamma)}(u|\tau,v), \lambda(y))$. Let $l(\psi, \alpha; \gamma)$ be the log-likelihood function from a single observed statistic and let $l_\alpha$ be the derivative of $l(\psi, \alpha; \gamma)$

with respect to $\alpha$ and $l_\psi$ be the differential operator of $l(\psi, \alpha; \gamma)$ with respect to $\psi$. According to the proof of Theorem 3.4, there exists a function $h(u, y, v) = (h_1(u, y, v), h_2(u, v), h_3(y))$ solving the equation $l_\psi^* l_\psi[h] = l_\psi^* l_\alpha$, where $l_\psi^*$ is the dual operator of $l_\psi$. Moreover, $\sqrt{n}(\hat{\alpha}_n - \alpha_0)$ is shown to have the asymptotic variance

$$(3.1) \quad E[\Sigma(\psi_0, \alpha_0, \gamma^*)^{-1}\Omega(O; \psi_0, \alpha_0, \gamma_0) + \omega(\psi_0, \alpha_0, \gamma_0) S(O; \gamma^*)]^{\otimes 2}.$$

Here, $\Sigma(\psi, \alpha, \gamma) = -\mathbf{P}[l_{\psi\alpha}[h] + l_{\alpha\alpha}]$ is the efficient information matrix for $\alpha$ for fixed $\gamma$, $\Omega(O; \psi, \alpha, \gamma) = l_\alpha - l_\psi[h]$ is the efficient score function for $\alpha$ for fixed $\gamma$, $\omega(\psi, \alpha, \gamma) = -\{\mathbf{P}[l_{\psi\alpha}[h] + l_{\alpha\alpha}]\}^{-1}\mathbf{P}[\nabla_\gamma(l_\psi[h] + l_\alpha)]$, and $S(O; \gamma^*)$ is the influence function of $\hat{\gamma}_n$.

To estimate (3.1), we wish to estimate each of the four terms including $\Sigma(\psi_0, \alpha_0, \gamma^*)$, $\Omega(O; \psi_0, \alpha_0, \gamma_0)$, $\omega(\psi_0, \alpha_0, \gamma_0)$ and $S(O; \gamma^*)$. At first, we define a pseudo-profile likelihood function as $pl_n(\alpha, \gamma) = n^{-1} \sum_{i=1}^n l_i(\hat{\psi}(\alpha, \gamma), \alpha; \gamma)$, where $l_i(\cdot)$ is the value of $l(\cdot)$ at the $i$th observation and $\hat{\psi}(\alpha, \gamma)$ is the argument of $\psi$ in the maximization of Section 2 when $\alpha$ and $\gamma$ are fixed. Then each of the four terms in (3.1) can be estimated using the following approach.

First, since $\Sigma(\psi_0, \alpha_0, \gamma^*)$ is the semiparametric efficiency information for $\alpha$ in the likelihood function of $(Y, R, U^*, V)$ when assuming $\gamma^*$ is known, according to Murphy and van der Vaart (2000), we can estimate it by

$$\hat{\Sigma}_n = -\frac{pl_n(\hat{\alpha}_n + \varepsilon_n, \hat{\gamma}_n) - 2pl_n(\hat{\alpha}_n, \hat{\gamma}_n) + pl_n(\hat{\alpha}_n - \varepsilon_n, \hat{\gamma}_n)}{\varepsilon_n^2}$$

where $\varepsilon_n$ is a constant of order $n^{-1/2}$.

Next, since $\hat{\psi}(\alpha, \hat{\gamma}_n)$ maximizes $\mathbf{P}_n l(\alpha, \psi; \hat{\gamma}_n)$, it holds that

$$\mathbf{P}_n[l_\psi(\hat{\psi}(\alpha, \hat{\gamma}_n), \alpha; \hat{\gamma}_n)[\tilde{h}]] = 0$$

for any tangent function $\tilde{h}$ of $\psi$. We differentiate the above equation with respect to $\alpha$, then evaluate it at $\hat{\alpha}_n$. This gives

$$\mathbf{P}_n l_{\alpha\psi}(\hat{\psi}(\hat{\alpha}_n, \hat{\gamma}_n), \hat{\alpha}_n; \hat{\gamma}_n)[\tilde{h}] = \mathbf{P}_n l_{\psi\psi}(\hat{\psi}(\hat{\alpha}_n, \hat{\gamma}_n), \hat{\alpha}_n; \hat{\gamma}_n)[\nabla_\alpha \hat{\psi}(\hat{\alpha}_n, \hat{\gamma}_n), \tilde{h}].$$

When $n$ goes to infinity, this equation approximates the equation which $h(u, y, v)$ solves. Thus, we expect that $\nabla_\alpha \hat{\psi}(\hat{\alpha}_n, \hat{\gamma}_n) \approx h(u, y, v)$. As a result, $\Omega(O_i; \psi_0, \alpha_0, \gamma^*) \approx \nabla_\alpha l_i(\hat{\psi}(\hat{\alpha}_n, \hat{\gamma}_n), \hat{\alpha}_n; \hat{\gamma}_n)$, while the latter can be evaluated using the numerical difference $\varepsilon_n^{-1}\{l_i(\hat{\psi}(\hat{\alpha}_n + \varepsilon_n, \hat{\gamma}_n), \hat{\alpha}_n + \varepsilon_n; \hat{\gamma}_n) - l_i(\hat{\psi}(\hat{\alpha}_n, \hat{\gamma}_n), \hat{\alpha}_n; \hat{\gamma}_n)\}$.

Third, we define $\hat{\alpha}(\gamma)$ as the estimate of $\alpha$ maximizing $pl_n(\alpha, \gamma)$ when $\gamma$ is held fixed. Using the argument similar to that in Zeng (2004), we can estimate $\omega(\psi_0, \alpha_0, \gamma^*)$ by a vector $\hat{\omega}_n$ with its $j$th element equal to $\tilde{\varepsilon}_n^{-1}(\hat{\alpha}(\hat{\gamma}_n + \tilde{\varepsilon}_n e_j) - \hat{\alpha}_n)$ for the $j$th canonical base $e_j$ and $\tilde{\varepsilon}_n$ satisfying $\tilde{\varepsilon}_n = o(1)$ and $\sqrt{n}\tilde{\varepsilon}_n \to \infty$.

Finally, $S(O; \gamma^*)$ can be estimated by $\hat{S}_n(O; \hat{\gamma}_n)$ using an explicit expression given in Zeng (2004).

Hence, the expression in (3.1) can be estimated by

$$\frac{1}{n} \sum_{i=1}^{n} \left[ \hat{\Sigma}_n^{-1} \frac{l_i(\hat{\psi}(\hat{\alpha}_n + \varepsilon_n), \hat{\alpha}_n + \varepsilon_n; \hat{\gamma}_n) - l_i(\hat{\psi}, \hat{\alpha}_n; \hat{\gamma}_n)}{\varepsilon_n} + \hat{\omega}_n \hat{S}_n(O_i, \hat{\gamma}_n) \right]^2.$$

**4. Simulation study.** A simulation study is conducted to illustrate our approach. In the simulation, for convenience of computation $V$ is chosen to be a binary variable with equal probabilities. Conditional on $V$, the lifetime $T$ is generated from a proportional hazards regression model with hazard rate $3t \exp\{V\}$. One surrogate variable $X_1$ is generated from the model $X_1 = \beta_0 T + 0.5\theta$, where $\theta$ is uniformly distributed in $(-0.5, 0.5)$ and $\beta_0$ may take different values in the simulation study. The study end time, $\tau$, is chosen to be 1. Additionally, we generate another irrelevant covariate $X_2$ from the uniform distribution in $[0, 1]$ and generate the right-censoring time from a proportional hazards model with hazard rate $4 \exp\{2X_1 - 4X_2 - 0.1V\}$. In other words, the simulation imitates the situation in which lifetime and censoring time are dependent and their dependence is explained by treatment $V$, a surrogate variable $X_1$ and a censorship related variable $X_2$.

According to our approach, the estimation of $\alpha$ is obtained by maximizing a pseudo-likelihood function over a sieve space, which is constructed similar to Section 2.3, with the choice $K_n = 5$ and $m = 3$ (other choices of $K_n$ and $m$ have little effect on the results, but large $K_n$ significantly increases computation time). Since $V$ is binary, for either value of $V$, $\eta_1(U, Y, V)$ is given as a linear combination of $N_{i_1}^m(U) N_{i_1}^m(Y)$ and $\eta_2(U, V)$ is given as a linear combination of $N_{i_1}^m(U)$. To prevent the parameters in the maximization from being unbounded, a penalty function, equal to $10^{-3}$ times the sum of squares of the spline coefficients, is subtracted from the pseudo-likelihood function. In the optimization, searching for the maximum starts from the initial values that $\alpha = 1$ and all the spline coefficients are zero. Our simulations show that the optimum search usually converges within 10 iterations when either the search-move step or the norm of the search direction is small enough.

The asymptotic variance of $\hat{\alpha}_n$ is estimated using the approach described in Section 3.3. Particularly, we choose $\varepsilon_n = n^{-1/2}, 3n^{-1/2}, 6n^{-1/2}$ and $\tilde{\varepsilon}_n = n^{-1/3}, 5n^{-1/3}$ in evaluating $\hat{\Sigma}_n^{-1}$ and $\hat{\omega}_n$. The results indicate that the estimates of the variance are pretty robust to these choices. Thus, only the results from $\varepsilon_n = n^{-1/2}$ and $\tilde{\varepsilon}_n = n^{-1/3}$ are reported here.

We choose $\beta_0 = 0$ or $\beta_0 = 1.5$ in the simulation. When $\beta_0 = 0$, the working model for $T$ is correct and the theoretical censoring rate is 18%; when $\beta_0 = 1.5$, the working model for $T$ is misspecified and the theoretical censoring rate becomes 36%. Table 1 summarizes the results from 500 repetitions

12                                D. ZENG

with sample size $n = 200$ for these two choices. In the table, the first column gives the true value of the parameter $\beta_0$. The second column gives the working models used in the estimation (e.g., $T|V$ means that the working model for $T$ is a proportional hazards model with $V$ as independent variable) and the superscript star in the column list indicates that the indexed working model is misspecified. In the third column, we report the naive estimates of $\alpha$ by regressing $T$ on $V$ directly. The remaining columns in turn report the average estimates of $\hat{\alpha}_n$, the standard errors of all the estimates, the median values of the estimated standard errors for $\hat{\alpha}_n$ and the coverage proportion of 95% confidence intervals based on the normal distribution approximation. Additionally, Figure 1 plots the histograms of $\hat{\alpha}_n$ from the simulations.

The simulation results indicate that when either working model is correct, the estimates produce small bias and moreover, our variance estimation approach gives fairly accurate estimates and valid coverage probabilities. Specifically, when $T$ is not fully predicted by $V$ and the working model for the censorship is correct, our estimate has smaller bias than the naive estimate. The simulation also shows that using the correct working model for $T$ may give a more efficient estimate. The amount of bias in $\hat{\alpha}_n$ observed in Table 1 can be due to the small sample size and the small $K_n$, as well as the imprecise evaluation of the integral in the likelihood function.

**5. Discussion.** For right-censored data, when the dependence between lifetime and censoring time is explained by many covariates, we utilize two working models to condense this high-dimensional information and thus derive the estimator of the treatment effect by maximizing some pseudo-likelihood function. We have shown that the estimator is consistent and asymptotically normal when either working model is correct.

For simplicity, the working model for $T$ given $W$ given in Section 2 is assumed to be the same as $T$ given $V$. This may seem very restrictive. However, in practice any semiparametric model can be adopted as the working model for $T$ given $W$. For example, suppose that we use a semiparametric model for $T$ given $W$ as follows: $f_{T|W}(y|w) = p(y, \beta^T w)$; then the condensed information will include $(U_1 = \beta^T W, U_2 = \gamma^T W, V)$. Hence, the estimator of $\alpha$ can be derived by maximizing

$$\mathbf{P}_n\{R\log[e^{\alpha V}\lambda(Y)\exp\{-\Lambda(Y)e^{\alpha V}\}f_{U_1,U_2}(\hat{U}_1,\hat{U}_2|Y,V)]\}$$

$$+ \mathbf{P}_n\bigg\{(1-R)\log\bigg[\int_Y^{\tau-} e^{\alpha V}\lambda(s)\exp\{-\Lambda(s)e^{\alpha V}\}f_{U_1,U_2}(\hat{U}_1,\hat{U}_2|s,V)\,ds$$

$$+ \exp\{-\Lambda(\tau)e^{\alpha V}\}f_{U_1,U_2}(\hat{U}_1,\hat{U}_2|T\geq\tau,V)\bigg]\bigg\}$$

over a sieve space of the parameters $(\alpha, f_{U_1,U_2}(u_1,u_2|T=y,v)$ $f_{U_1,U_2}(u_1,u_2|T\geq \tau,v), \lambda(y))$, where $\hat{U}_1 = \hat{\beta}_n^T W$ and $\hat{U}_2 = \hat{\gamma}_n^T W$ for some estimators $\hat{\beta}_n, \hat{\gamma}_n$. The slight difference from the previous context is that B-splines in the sieve space are constructed on a four-dimensional space. Consequently, under some regular conditions, one of the following two conclusions is expected to be true: if the semiparametric working model for $T$ given $W$ does not satisfy the constraint that $h_{T|V}(y|v) = \lambda(y)e^{\alpha^T v}$, that is, the working model is misspecified, the consistency of $\hat{\alpha}_n$ holds if the working model for $C$ given $W$ is correct;

TABLE 1
*Simulation results from* 500 *repetitions with sample size* 200

| $\beta_0$ | Working models | Naive est. | $\overline{\hat{\alpha}_n}$ | $\text{se}(\hat{\alpha}_n)$ | $\text{med}(\widehat{\text{se}})$ | 95% CI |
|---|---|---|---|---|---|---|
| 0 | $(T|V), (C|X_1, X_2, V)$ | 1.004 | 0.974 | 0.169 | 0.172 | 0.956 |
|   | $(T|V), (C|X_2, V)^*$ | 1.004 | 0.975 | 0.169 | 0.172 | 0.954 |
| 1.5 | $(T|V)^*, (C|X_1, X_2, V)$ | 0.835 | 0.915 | 0.189 | 0.234 | 0.976 |
|   | $(T|V)^*, (C|X_2, V)^*$ | 0.835 | 0.802 | 0.186 | 0.208 | 0.868 |

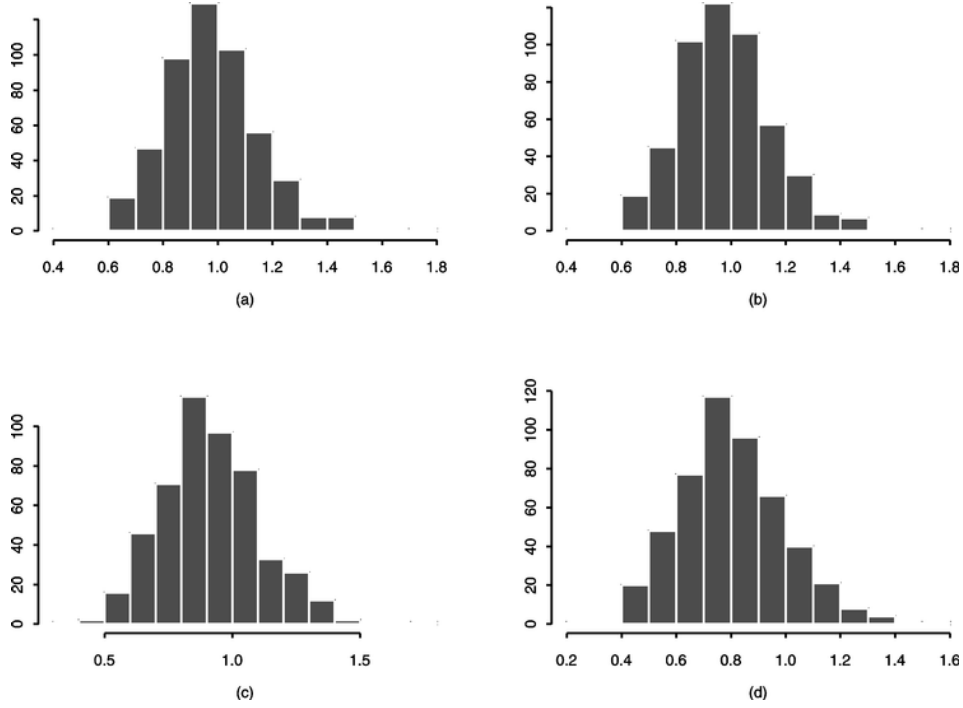

FIG. 1. *Histograms of $\hat{\alpha}_n$ from 500 repetitions:* (a) *both working models for $T$ and $C$ are correct,* (b) *the working model for $T$ is correct but the working model for $C$ is misspecified,* (c) *the working model for $T$ is misspecified but the working model for $C$ is correct,* (d) *both working models are misspecified.*

on the contrary, if the working model for $T$ given $W$ satisfies the constraint, the double robustness of $\hat{\alpha}_n$ given in Theorem 3.4 holds as well. However, it is often difficult to specify a correct working model for $T$ given $W$ satisfying the constraint (1.1) except in the simplest situation that $T$ depends on $W$ only via $V$, which has been used in this paper.

Our approach can be easily extended to the situation when $V$ is multidimensional and possibly discrete. If $V$ is multidimensional, the sieve space needs to be constructed on a real space of all of $U$, $Y$ and the multidimensional $V$. However, if $V$ is discrete, the sieve space only needs to be constructed on a real space $U$ and $Y$ for each category of $V$. The latter has already been implemented in the simulation study.

We acknowledge that our approaches are not easily generalized to the situation with a time-dependent component in $X$, since when $X$ contains time-dependent covariates the condensed information using working models still depends on time, so its dimensionality is not reduced essentially. Further investigation is being conducted to solve this problem.

**6. Proofs.** For convenience of writing, we assume $\tau = 1$ and denote $\mathbf{G}_n$ as the empirical process $\sqrt{n}(\mathbf{P}_n - \mathbf{P})$.

PROOF OF THEOREM 3.1. The whole proof can be divided into three steps: first, we construct some functions in the sieve space which approximate the true parameters; then by using empirical process theory, we obtain one key inequality; finally, this inequality is used to obtain the consistency.

*Step* 1. We construct some functions in $\mathcal{S}_n(m, K_n, M_n)$ to approximate the true parameters. To do that, we need the following general result. From the properties of B-spline functions [cf. Schumaker (1981)], we can define a linear operator $\mathcal{Q}_p$ mapping $W^{k,\infty}([0,1]^p)$ to the sieve space; that is, for any $g \in W^{k,\infty}([0,1]^p)$,

$$\mathcal{Q}_p[g] = \sum_{i_1,\ldots,i_p=1}^{m+K_n} \Gamma_{i_1,\ldots,i_p}[g] N_{i_1}^m(x_1) \ldots N_{i_p}^m(x_p),$$

where $\Gamma_{i_1,\ldots,i_p}$ are the linear functionals in $L_\infty([0,1]^p)$. Moreover,

$$\sum_{i_1,\ldots,i_p}^{m+K_n} |\Gamma_{i_1,\ldots,i_p}[g]| \leq (2m+1)^p 9^{p(m-1)} \|g\|_{L_\infty([0,1]^p)},$$

and according to Schumaker [(1981), Theorem 12.7],

$$\|\mathcal{Q}_p[g] - g\|_{L_\infty([0,1]^p)} \leq \frac{C(m)}{K_n^k} \|g\|_{W^{k,\infty}([0,1]^p)}.$$

Thus, we define $\eta_{1n}(u,y,v) = \mathcal{Q}_3[\log f_{U^*}] - \mathcal{Q}_3[\log f_{U^*}]|_{u=0}$, $\eta_{2n}(u,v) = \mathcal{Q}_2[\log g_{U^*}] - \mathcal{Q}_2[\log g_{U^*}]|_{u=0}$ and $\xi_n(y) = \mathcal{Q}_1[\log \lambda_0]$. Correspondingly, we obtain

$$f_n(u|y,v) = \frac{\exp\{\eta_{1n}(u,y,v)\}}{\int_0^1 \exp\{\eta_{1n}(u,y,v)\}\,du}, \qquad g_n(u|\tau,v) = \frac{\exp\{\eta_{2n}(u,v)\}}{\int_0^1 \exp\{\eta_{2n}(u,v)\}\,du},$$

and $\lambda_n(y) = \exp\{\xi_n(y)\}$. As a result of the fact that $\sum_{i_1=1}^{m+K_n} N_{i_1}^m(u) = 1$, $(\alpha_0, \xi_n(y), \eta_{1n}(u,y,v), \eta_{2n}(u,v))$ is in the sieve space $\mathcal{S}_n(m, K_n, M_n)$ and moreover,

$$\|f_n - f_{U^*}\|_{L_\infty([0,1]^3)} \leq O(1)\|\log f_{U^*} - \mathcal{Q}_3[\log f_{U^*}]\|_{L_\infty([0,1]^3)} \leq O\left(\frac{1}{K_n^k}\right)$$

and the same bound holds for $\|g_n - g_{U^*}\|_{L_\infty([0,1]^2)}$ and $\|\lambda_n - \lambda_0\|_{L_\infty([0,1])}$.

*Step* 2. We obtain a key inequality using empirical process theory. To simplify the notation, for any functions $f_1(u,y,v)$, $f_2(u,v)$ and $f_3(y)$, we denote $G(r, f_1, f_2, f_3, \alpha; \gamma)$ as the likelihood function from one single observation with parameters $(\alpha, f_3, f_1, f_2)$. Since $(\hat{\alpha}_n, \hat{\lambda}_n, \hat{f}_n, \hat{g}_n)$ maximizes $\mathbf{P}_n[\log G(R, f_1, f_2, f_3, \alpha; \hat{\gamma}_n)]$ over the sieve space, it follows that

$$\mathbf{P}_n[\log G(R, \hat{f}_n, \hat{g}_n, \hat{\lambda}_n, \hat{\alpha}_n; \hat{\gamma}_n)] \geq \mathbf{P}_n[\log G(R, f_n, g_n, \lambda_n, \alpha_0; \hat{\gamma}_n)].$$

Equivalently,

$$
\begin{aligned}
n^{-1/2}\mathbf{G}_n&\left[\log \frac{G(R, \hat{f}_n, \hat{g}_n, \hat{\lambda}_n, \hat{\alpha}_n; \hat{\gamma}_n)}{G(R, f_n, g_n, \lambda_n, \alpha_0; \hat{\gamma}_n)}\right] \\
&\geq \mathbf{P}\left[\log \frac{G(R, f_n, g_n, \lambda_n, \alpha_0; \hat{\gamma}_n)}{G(R, f_{U^*}, g_{U^*}, \lambda_0, \alpha_0; \gamma^*)}\right] \\
&\quad + \mathbf{P}\left[\log \frac{G(R, f_{U^*}, g_{U^*}, \lambda_0, \alpha_0; \gamma^*)}{G(R, \hat{f}_n, \hat{g}_n, \hat{\lambda}_n, \hat{\alpha}_n; \hat{\gamma}_n)}\right],
\end{aligned}
\tag{6.1}
$$

where we recall that $f_{U^*}$ and $g_{U^*}$ are the conditional densities of $U^*$ given $(T, V)$ and $(T \geq \tau, V)$, respectively.

We want to bound the left-hand side of (6.1) using empirical process theory. For this purpose, we consider a class of functions $\mathcal{L}_n$ defined by

$$\left\{\log \frac{G(r, \tilde{f}_n, \tilde{g}_n, \tilde{\lambda}_n, \alpha; \hat{\gamma}_n)}{G(r, f_n, g_n, \lambda_n, \alpha_0; \hat{\gamma}_n)} : \tilde{\lambda}_n(y) = e^{\tilde{\xi}(y)},\right.$$

$$\tilde{f}_n(u|y,v) = \frac{\exp\{\tilde{\eta}_1(u,y,v)\}}{\int_0^1 \exp\{\tilde{\eta}_1(\tilde{u},y,v)\}\,d\tilde{u}},$$

$$\tilde{g}_n(u|\tau,v) = \frac{\exp\{\tilde{\eta}_2(u,v)\}}{\int_0^1 \exp\{\tilde{\eta}_2(\tilde{u},v)\}\,d\tilde{u}},$$

$$\left.(\alpha, \tilde{\xi}, \tilde{\eta}_1, \tilde{\eta}_2) \in \mathcal{S}_n(m, K_n, M_n)\right\}.$$

Since $\|N_i^m(.)\|_{L_\infty([0,1])} = 1$, any function of $\tilde{f}_n, \tilde{g}_n, \tilde{\lambda}_n$ given in $\mathcal{L}_n$ is bounded by $O(e^{2M_n})$. By Assumptions A1 and A2, $G(r, f_n, g_n, \lambda_n, \alpha_0; \hat{\gamma}_n)$ is bounded away from 0. Hence, the class $\mathcal{L}_n$ has an upper bound $O_p(M_n)$. Moreover, this class can be regarded as the class of functions indexed by $\alpha$, $\{\tilde{\eta}^1_{i_1,i_2,i_3}\}_{i_1,i_2,i_3=1}^{m+K_n}$, $\{\tilde{\eta}^2_{i_1,i_2}\}_{i_1,i_2=1}^{m+K_n}$ and $\{\tilde{\xi}_i\}_{i=1}^{m+K_n}$, which are the respective B-spline coefficients of $\tilde{\eta}_1, \tilde{\eta}_2$ and $\tilde{\xi}$ in $\mathcal{S}_n(m, K_n, M_n)$. Tedious checking indicates that the function in $\mathcal{L}_n$ is Lipschitz continuous with respect to all these parameters and the Lipschitz constant is bounded by $O_p(e^{6M_n})$. In addition, since $|\tilde{\eta}^1_{i_1,i_2,i_3}|, |\tilde{\eta}^2_{i_1,i_2}|$ and $|\tilde{\xi}_i|$ are bounded by $M_n$ and $|\alpha|$ is bounded by $M$, they lie in a hypercube of a real space $R^{N_n+1}$ where $N_n = (m+K_n)^3 + (m+K_n)^2 + m + K_n$. Therefore, for any $\varepsilon > 0$, if we partition this hypercube into subcubes with scale length $\varepsilon$, the total number of subcubes is at most $O((M_n/\varepsilon)^{N_n})$. According to the Lipschitz property of the functions in $\mathcal{L}_n$, the $L_\infty$-distance between any two functions of $\mathcal{L}_n$ with respective indexes in the same subcube is no more than $O_p(e^{6M_n})N_n\varepsilon$. Consequently, we obtain that the bracketing number for $\mathcal{L}_n$ satisfies $N_{[\cdot]}(O_p(e^{6M_n})N_n\varepsilon, \mathcal{L}_n, L_\infty) \leq O(1)(M_n/\varepsilon)^{N_n}$. According to van der Vaart [(1998), Theorem 19.35], in probability we have

$$\sqrt{n} E_P^* \|\mathbf{P}_n - \mathbf{P}\|_{\mathcal{L}_n} \leq O_p(1) \int_0^{O(M_n)} \sqrt{\log\left(\frac{2M_n e^{6M_n}(m+K_n)^3}{\varepsilon}\right)^{2(m+K_n)^3}} d\varepsilon$$

$$\leq O_p(1) K_n^{3/2} (\log K_n) M_n^2.$$

Thus, the left-hand side of inequality (6.1) is bounded by $O_p(M_n^2 K_n^{3/2} \log K_n/\sqrt{n})$ from above.

We denote the two terms in the right-hand side of (6.1) as $(I)$ and $(II)$ and wish to bound them from below. Since the functional $G(\cdot)$ is Lipschitz continuous with each component, we have that

$$(I) \geq -O_p(1)\{\|f_n - f_{U^*}\|_{L_\infty} + \|g_n - g_{U^*}\|_{L_\infty} + \|\lambda_n - \lambda_0\|_{L_\infty} + |\hat{\gamma}_n - \gamma^*|\}$$

$$\geq -O_p(1)\left(\frac{1}{K_n^k} + \frac{1}{\sqrt{n}}\right).$$

On the other hand, by Schumaker [(1981), Theorem 4.22] $|dN_{i_1}^m(u)/du| \leq O(K_n)$. We can easily verify that

$$|G(R, \hat{f}_n, \hat{g}_n, \hat{\xi}, \hat{\alpha}_n; \hat{\gamma}_n) - G(R, \hat{f}_n, \hat{g}_n, \hat{\xi}, \hat{\alpha}_n; \gamma^*)| \leq O(e^{2M_n} M_n K_n)|\hat{\gamma}_n - \gamma^*|.$$

Therefore,

$$(II) \geq -O_p(e^{2M_n}) M_n K_n |\hat{\gamma}_n - \gamma^*| + \mathbf{P}\left[\log \frac{G(R, f_{U^*}, g_{U^*}, \lambda_0, \alpha_0; \gamma^*)}{G(R, \hat{f}_n, \hat{g}_n, \hat{\lambda}_n, \hat{\alpha}_n; \gamma^*)}\right].$$

However, the last term in the above is the Kulback–Leibler information. We linearize the last term. The first-order term in the expansion vanishes while the second-order term in the expansion is bounded from below by

$$O(e^{-3M_n})\|G(R, f_{U^*}, g_{U^*}, \lambda_0, \alpha_0; \gamma^*) - G(R, \hat{f}_n, \hat{g}_n, \hat{\lambda}_n, \hat{\alpha}_n; \gamma^*)\|^2_{L_2(P)}.$$

Combining the above results and noting that the probability measure $P$ is equivalent to the product measure of the Lebsgue measure in $[0,1]^3$ and the counting measure for $\{0,1\}$, we obtain that for $r = 0, 1$,

$$
\begin{aligned}
O_p(1)&\left(\frac{e^{5M_n}M_n K_n}{\sqrt{n}} + \frac{e^{3M_n}}{K_n^k} + \frac{e^{3M_n}M_n^2 K_n^{3/2}\log K_n}{\sqrt{n}}\right) \\
&\geq \int_{[0,1]^3} [G(r, f_{U^*}, g_{U^*}, \lambda_0, \alpha_0; \gamma^*) \\
&\qquad\qquad - G(r, \hat{f}_n, \hat{g}_n, \hat{\lambda}_n, \hat{\alpha}_n; \gamma^*)]^2 \, du\, dy\, dv.
\end{aligned}
$$
(6.2)

*Step* 3. We obtain the $L_2$-convergence of the estimators. Suppose we select $K_n$ and $M_n$ such that they satisfy Assumption A5. Equation (6.2) implies that this upper bound holds for the square $L_2$-distance between $\int_0^s \int_0^1 G(1, f_{U^*}, g_{U^*}, \lambda_0, \alpha_0; \gamma^*)\, du\, dy$ and $\int_0^s \int_0^1 G(1, \hat{f}_n, \hat{g}_n, \hat{\lambda}_n, \hat{\alpha}_n; \gamma^*)\, du\, dy$ for any $s \in [0,1]$. After simplification, we obtain that

$$
\begin{aligned}
\int &[\exp\{-e^{\hat{\alpha}_n v}\hat{\Lambda}_n(s)\} - \exp\{-e^{\alpha_0 v}\Lambda_0(s)\}]^2\, dv \\
&\leq O_p(1)\left(\frac{e^{3M_n}}{K_n^k} + \frac{e^{6M_n}K_n^{3/2}\log K_n}{\sqrt{n}}\right).
\end{aligned}
$$
(6.3)

By choosing a subsequence, we suppose $\hat{\alpha}_n \to \alpha^*$. From the above inequality and Assumption A1, $\alpha^* = \alpha_0$ and $\hat{\Lambda}_n(y)$ converges pointwise to $\Lambda_0(y)$ for $y \in [0,1]$. Furthermore, since $\Lambda_0$ is continuous, $\|\hat{\Lambda}(y) - \Lambda_0(y)\|_{L_\infty([0,1])} \to 0$. This completes the proof of Theorem 3.1. $\square$

PROOF OF THEOREM 3.2. From the last inequality and Assumption A1, we immediately obtain that

$$|\hat{\alpha}_n - \alpha_0|^2 \leq O_p(1)\left(\frac{e^{3M_n}}{K_n^k} + \frac{e^{6M_n}K_n^{3/2}\log K_n}{\sqrt{n}}\right).$$

After repeating using (6.2) for $R = 1$ and $R = 0$, we can further obtain that the same bound holds for $\|\hat{\lambda}_n - \lambda_0\|^2_{L_2(P)}$, $\|\hat{f}_n - f_{U^*}\|^2_{L_2(P)}$ and $\|\hat{g}_n - g_{U^*}\|^2_{L_2(P)}$.

On the other hand, from Schumaker [(1981), Theorem 4.22] and Assumption A3, we have that

$$\|\nabla_u^{k_1}\nabla_y^{k_2}\nabla_v^{k_3}\hat{\eta}_{1n}(u,y,v)\|_{L_\infty(P)} \leq CK_n^k \sum_{i_1,i_2,i_3=1}^{m+K_n} |\hat{\eta}_{i_1,i_2,i_3}| \leq O(M_n K_n^k),$$

where $k_1+k_2+k_3=k$. Thus, $\|\nabla_u^{k_1}\nabla_y^{k_2}\nabla_v^{k_3}\hat{f}_n(u|y,v)\|_{L_\infty(P)} \leq Ce^{(k+1)M_n}M_n K_n^k$. According to the Sobolev interpolation inequality [cf. Adams (1975)], we obtain that

$$\|\nabla(\hat{f}_n - f_{U^*})\|_{L_\infty(P)} \leq Ce^{(k+2)M_n\tau_1} K_n^{k\tau_1}\left(\frac{e^{3M_n}}{K_n^k} + \frac{e^{6M_n}K_n^{3/2}\log K_n}{\sqrt{n}}\right)^{(1-\tau_1)/2},$$

where $\tau_1 = \frac{5/2}{k}$. By the choice of $K_n$ and $M_n$ in Assumption A5, $\|\nabla(\hat{f}_n - f_{U^*})\|_{L_\infty(P)}$ converges to zero. Similarly, this is true for $\hat{g}_n$ and $\hat{\lambda}_n$. Thus, Theorem 3.2 holds. $\square$

PROOF OF THEOREM 3.3. Using the results from Theorems 3.1 and 3.2, redo the proof of Theorem 3.1. We define $\mathcal{L}_n$ as a class as before, but the functions in $\mathcal{L}_n$ are indexed by $(\alpha, \xi, f_{U(\gamma)}, g_{U(\gamma)})$, which belongs to a bounded set in $R \times \{W^{1,\infty}(P)\}^3$. Thus, $\mathcal{L}_n$ has a bounded covering function and the integration of the entropy for the class $\mathcal{L}_n$ is finite. Moreover, the function in the left-hand side of (6.1) converges to zero uniformly. Thus, we can apply Theorem 2.11.23 of van der Vaart and Wellner (1996), to obtain that the left-hand side of inequality (6.1) is bounded by $o_p(1/\sqrt{n})$. For the right-hand side of (6.1), we still perform Taylor expansion at the true parameters. Since each parameter is in a small neighborhood of the true parameters, the right-hand side of (6.1) is bounded from below by

$$-O_p\{|\hat{\gamma}_n - \gamma^*|^2 + \|f_n - f_{U^*}\|_{L_2(P)}^2 + \|g_n - g_{U^*}\|_{L_2(P)}^2 + \|\lambda_n - \lambda_0\|_{L_2(P)}^2\}$$
$$+ O_p(1)\|G(R, f_{U^*}, g_{U^*}, \lambda_0, \alpha_0; \gamma^*) - G(R, \hat{f}_n, \hat{g}_n, \hat{\lambda}_n, \hat{\alpha}_n; \gamma^*)\|_{L_2(P)}^2.$$

Recall the construction of $f_n, g_n$ and $\xi_n$ in the first step of proving Theorem 3.1; we obtain that

$$\frac{o_p(1)}{\sqrt{n}} + \frac{O_p(1)}{K_n^{2k}}$$
$$\geq \|G(R, f_{U^*}, g_{U^*}, \lambda_0, \alpha_0; \gamma^*) - G(R, \hat{f}_n, \hat{g}_n, \hat{\lambda}_n, \hat{\alpha}_n; \gamma^*)\|_{L_2(P)}^2.$$

The results of Theorem 3.3 thus follow from the same arguments as in the proof of Theorem 3.2. $\square$

PROOF OF THEOREM 3.4. We will write $\sqrt{n}(\hat{\alpha}_n - \alpha_0)$ as a linear functional of the empirical process $\mathbf{G}_n$. The whole proof can be divided into the following five steps.

*Step* 1. We define a pseudo least favorable direction for $\alpha_0$ when $\gamma^*$ is known. The nuisance parameters for $\alpha$ are $(f_{U^*}, g_{U^*}, \lambda)$ and are denoted as $\psi$. The tangent space for $\psi$ is thus given by

$$H = \left\{ h(u,y,v) = (h_1(u,y,v), h_2(u,v), h_3(y)) : \int_0^1 h_1(u,y,v)\, du = 0, \right.$$
$$\left. \int_0^1 h_2(u,v)\, du = 0,\ h(u,y,v) \in L_2([0,1]^3) \right\}.$$

Let $l(\psi, \alpha; \gamma^*)$ be as defined in Section 3.3. Then a pseudo least favorable direction for $\alpha_0$ is defined as a tangent function $h(u,y,v) \in H$ for $\psi$ that satisfies

$$l_\psi^*(\psi_0, \alpha_0; \gamma^*) l_\psi(\psi_0, \alpha_0; \gamma^*)[h] = l_\psi^*(\psi_0, \alpha_0; \gamma^*) l_\alpha(\psi_0, \alpha_0; \gamma^*) \qquad \text{a.s.},$$

where $l_\psi(\psi_0, \alpha_0; \gamma^*)[h]$ is the derivative of $l(\cdot)$ with respect to $\psi$ along the direction $h$, $l_\psi^*(\psi_0, \alpha_0; \gamma^*)$ is the adjoint operator of $l_\psi(\psi_0, \alpha_0; \gamma^*)$ in the Hilbert space $L_2(P)$ and $l_\alpha(\psi_0, \alpha_0; \gamma^*)$ is the derivative of $l$ respective to $\alpha$.

*Step* 2. We prove the existence of the pseudo least favorable direction. We note that $H$ is a Hilbert space with $\langle h, h \rangle_H$ given by

$$\|h_1(u,y,v)\|^2_{L_2([0,1]^3)} + \|h_2(u,v)\|^2_{L_2([0,1]^2)} + \|h_3(y)\|^2_{L_2([0,1])}.$$

Then the following lemma holds.

LEMMA 6.1. *Under Assumptions* A1–A4, *there exists a unique* $h \in H$ *such that*

$$l_\psi^*(\psi_0, \alpha_0; \gamma^*) l_\psi(\psi_0, \alpha_0; \gamma^*)[h] = l_\psi^*(\psi_0, \alpha_0; \gamma^*) l_\alpha(\psi_0, \alpha_0; \gamma^*) \qquad \text{a.s.}$$

PROOF. Define $A$ as a linear operator from $L_2([0,1])$ to $L_2([0,1]^2)$, given by $A[h_3] = -e^{\alpha_0 v} \int_0^y h_3(s)\, ds + h_3(y)/\lambda_0(y)$. After some calculation and using the property $\int_0^1 h_1(u,y,v)\, du = 0$, we have

$$\|l_\psi(\psi_0, \alpha_0; \gamma^*)[h]\|^2_{L_2(P)}$$
$$\geq \left\| R\left[ A[h_3] + \frac{h_1(U^*, Y, V)}{f_{U^*}(U^*|Y,V)} \right] \right\|^2_{L_2(P)}$$
$$+ \left\| (1-R) I(Y=\tau) \frac{h_2(U^*, V)}{g_{U^*}(U^*|\tau, V)} \right\|^2_{L_2(P)}$$
$$\geq O(1) \left[ \int_{[0,1]^2} A[h_3]^2\, dy\, dv + \|h_1\|^2_{L_2([0,1]^3)} \right] + \|h_2\|^2_{L_2([0,1]^2)}.$$

Since $A$ is invertible, $\|A[h_3]\| \geq \|A^{-1}\|^{-1} \|h_3\|$. Thus, the last term is bounded from below by $O(1)\langle h, h \rangle_H$. By the Lax–Milgram theorem [Evans (1998)], the operator $l_\psi^*(\psi_0, \alpha_0; \gamma^*) l_\psi(\psi_0, \alpha_0; \gamma^*)$ is invertible. Lemma 6.1 is proved. □

*Step* 3. The proof for the smoothness of the least favorable direction is technical so we leave it to one of our technical reports, which is available from the author. There we show:

LEMMA 6.2.  *Under Assumptions* A1–A4, $h(u, y, v) \in W^{k,\infty}(P)$.

*Step* 4. We construct the projection of $h(u, y, v)$ on the tangent space of the sieve space. First, by simple computation, the tangent vectors $h_n(u, y, v)$ for the nuisance parameters at $\hat{\psi}_n = (\hat{f}_n(u|y, v), \hat{g}_n(u|\tau, v), \hat{\lambda}_n(y))$ have the form

$$\left( \hat{f}_n(u|y,v)\xi_1(u,y,v) - \hat{f}_n(u|y,v)\frac{\int_0^1 \exp\{\hat{\eta}_{1n}(u,y,v)\}\xi_1(u,y,v)\,du}{\int_0^1 \exp\{\hat{\eta}_{1n}(u,y,v)\}\,du}, \right.$$

$$\left. \hat{g}_n(u|\tau,v)\xi_2(u,v) - \hat{g}_n(u|\tau,v)\frac{\int_0^1 \exp\{\hat{\eta}_{2n}(u,v)\}\xi_2(u,v)\,du}{\int_0^1 \exp\{\hat{\eta}_{2n}(u,v)\}\,du}, \hat{\lambda}_n(y)\xi_3(y) \right),$$

where $\xi_1(u, y, v)$, $\xi_2(u, v)$ and $\xi_3(y)$ have the same forms as $\eta_1(u, y, v)$, $\eta_2(u, y, v)$ and $\xi(y)$ in the sieve space. Then, one good approximation to the pseudo least favorable direction is to choose $h_n(u, y, v) = (h_1^n, h_2^n, h_3^n)$ so that their corresponding $(\xi_1(u,y,v), \xi_2(u,v), \xi_3(y))$ satisfy $\xi_1(u,y,v) = \mathcal{Q}_3[h_1/f_{U^*}] - \mathcal{Q}_3[h_1/f_{U^*}]|_{u=0}$, $\xi_2(u,v) = \mathcal{Q}_2[h_2/g_{U^*}] - \mathcal{Q}_2[h_2/g_{U^*}]|_{u=0}$ and $\xi_3(y) = \mathcal{Q}_1[h_3/\lambda_0]$. Here the operator $\mathcal{Q}_p$ was defined in the proof of Theorem 3.1. Thus, the results in Theorem 3.3 and Lemma 6.2 imply that

$$\|h_n(U^*, Y, V) - h(U^*, Y, V)\|_{L_2(P)}^2 \le O\left(\frac{1}{K_n^{2k}}\right) + o_p\left(\frac{1}{\sqrt{n}}\right).$$

*Step* 5. We derive the empirical process for $\sqrt{n}(\hat{\alpha}_n - \alpha_0)$. Since $(\hat{\psi}_n, \hat{\alpha}_n)$ maximizes the log-likelihood in the sieve space, the score along the path $(\hat{\alpha} + \varepsilon, \hat{\psi} + \varepsilon h_n)$ is zero when $\varepsilon = 0$. Then it holds that

$$\mathbf{G}_n\{l_\psi(\hat{\psi}_n, \hat{\alpha}_n; \hat{\gamma}_n)[h_n] + l_\alpha(\hat{\psi}_n, \hat{\alpha}_n; \hat{\gamma}_n)\}$$
$$= -\sqrt{n}\mathbf{P}\{l_\psi(\hat{\psi}_n, \hat{\alpha}_n; \hat{\gamma}_n)[h_n] + l_\alpha(\hat{\psi}_n, \hat{\alpha}_n; \hat{\gamma}_n)\}.$$

For the left-hand side of the above equation, we apply Theorem 2.11.23 of van der Vaart and Wellner (1996). Note that the function in the left-hand side, indexed by both $(\hat{\psi}_n, h_n) \in W^{1,\infty}$ and $(\hat{\alpha}_n, \hat{\gamma}_n) \in [-M, M]^{d+2}$, belongs to a $P$-Donsker class. Moreover, we linearize the right-hand side at the true parameters and approximate $h_n$ by $h$. Since $\mathbf{P}\{l_{\psi\psi}(\psi_0, \alpha_0; \gamma^*)[\hat{\psi}_n - \psi_0, h] + l_{\alpha\psi}(\psi_0, \alpha_0; \gamma^*)[\hat{\psi}_n - \psi_0]\} = 0$, we obtain that

$$-\mathbf{P}\{l_{\psi\alpha}(\psi_0, \alpha_0; \gamma^*)[h] + l_{\alpha\alpha}(\psi_0, \alpha_0; \gamma^*)\}\sqrt{n}(\hat{\alpha}_n - \alpha_0)$$
$$= \mathbf{G}_n\{l_\psi(\psi_0, \alpha_0; \gamma^*)[h] + l_\alpha(\psi_0, \alpha_0; \gamma^*)\}$$

$$+ \mathbf{P}\{l_{\psi\gamma}(\psi_0, \alpha_0; \gamma^*)[h] + l_{\alpha\gamma}(\psi_0, \alpha_0; \gamma^*)\}\sqrt{n}(\hat{\gamma}_n - \gamma^*)$$
$$+ \sqrt{n}O_p(\|\hat{\psi}_n - \psi_0\|^2_{L_2(P)} + |\hat{\alpha} - \alpha_0|^2 + \|h_n - h\|^2_{L_2(P)} + |\hat{\gamma}_n - \gamma^*|^2).$$

The last term is $o_p(1)$ from Theorem 3.3 and Assumption A6. Hence, the asymptotic normality of $\sqrt{n}(\hat{\alpha}_n - \alpha_0)$ holds if we can prove the following lemma.

LEMMA 6.3. $-\mathbf{P}\{l_{\psi\alpha}(\psi_0, \alpha_0; \gamma^*)[h] + l_{\alpha\alpha}(\psi_0, \alpha_0; \gamma^*)\} > 0.$

PROOF. We note that

$$-\mathbf{P}\{l_{\psi\alpha}(\psi_0, \alpha_0; \gamma^*)[h] + l_{\alpha\alpha}(\psi_0, \alpha_0; \gamma^*)\}$$
$$= \mathbf{P}\{l_\alpha(\psi_0, \alpha_0; \gamma^*) + l_\psi(\psi_0, \alpha_0; \gamma^*)[h]\}^2 \geq 0.$$

Moreover, if $l_\alpha(\psi_0, \alpha_0; \gamma^*) + l_\psi(\psi_0, \alpha_0; \gamma^*)[h]$ is zero, then for $R=1$,

$$0 = \frac{h_1(U^*, Y, V)}{f_{U^*}(U^*|Y, V)} + \left\{\frac{h_3(Y)}{\lambda_0(Y)} - e^{\alpha_0 V}\int_0^Y h_3(s)\,ds\right\} - Ve^{-\alpha_0 V}\Lambda_0(Y) + V.$$

Let $Y = 0$. Multiply both sides by $f_{U^*}$, then integrate both sides over $U^*$ from $a(\gamma^*)$ to $b(\gamma^*)$. We have $V = -h_3(0)/\lambda_0(0)$. This is a contradiction. □

Furthermore, we obtain the influence function of $\hat{\alpha}_n$ to be

$$-\{\mathbf{P}[l_{\psi\alpha}(\psi_0, \alpha_0; \gamma^*)[h] + l_{\alpha\alpha}(\psi_0, \alpha_0; \gamma^*)]\}^{-1}$$
$$\times \{l_\psi(\psi_0, \alpha_0; \gamma^*)[h] + l_\alpha(\psi_0, \alpha_0; \gamma^*)$$
$$+ \mathbf{P}[l_{\psi\gamma}(\psi_0, \alpha_0; \gamma^*)[h] + l_{\alpha\gamma}(\psi_0, \alpha_0; \gamma^*)]S(Y, R, W; \gamma^*)\},$$

where $S(Y, R, W; \gamma^*)$ is the influence function of $\hat{\gamma}_n$. When both working models are correct, $l(\psi_0, \alpha_0; \gamma)$ is always the logarithm of the density for $(T \wedge C, R, U(\gamma), V)$ whatever value $\gamma$ takes. So $\mathbf{P}[l_\psi(\psi_0, \alpha_0; \gamma)[h] + l_\alpha(\psi_0, \alpha_0; \gamma)]$ is an expectation of a score function; thus, it is equal to 0. This implies $\mathbf{P}[l_{\psi\gamma}(\psi_0, \alpha_0; \gamma^*)[h] + l_{\alpha\gamma}(\psi_0, \alpha_0; \gamma^*)] = 0$. Hence, $\hat{\alpha}_n$ has an influence function equal to $-\{\mathbf{P}[l_{\psi\alpha}(\psi_0, \alpha_0; \gamma^*)[h] + l_{\alpha\alpha}(\psi_0, \alpha_0; \gamma^*)]\}^{-1}\{l_\psi(\psi_0, \alpha_0; \gamma^*)[h] + l_\alpha(\psi_0, \alpha_0; \gamma^*)\}$, which is exactly the efficient influence function for $\alpha$. Consequently, the asymptotic variance of $\sqrt{n}(\hat{\alpha}_n - \alpha_0)$ is equal to the semiparametric efficiency bound. □

**Acknowledgments.** This work is part of my Ph.D. dissertation advised by Professor Susan Murphy at the University of Michigan. I owe many thanks to her for numerous discussions and helpful comments. I also thank an Associate Editor and a referee for their valuable suggestions.

# REFERENCES

ADAMS, R. A. (1975). *Sobolev Spaces*. Academic Press, New York. MR450957

BICKEL, P. J., KLAASSEN, C. A. J., RITOV, Y. and WELLNER, J. A. (1993). *Efficient and Adaptive Estimation for Semiparameric Models*. Johns Hopkins Univ. Press. MR1245941

DEIMLING, K. (1985). *Nonlinear Functional Analysis*. Springer, Berlin. MR787404

EVANS, L. C. (1998). *Partial Differential Equations*. Amer. Math. Soc., Providence, RI. MR1625845

MURPHY, S. A. and VAN DER VAART, A. W. (2000). On profile likelihood (with discussion). *J. Amer. Statist. Assoc.* **95** 449–485. MR1803168

ROBINS, J. M., ROTNITZKY, A. and VAN DER LAAN, M. (2000). Comment on "On profile likelihood," by S. Murphy and A. W. van der Vaart. *J. Amer. Statist. Assoc.* **95** 477–482.

ROBINS, J. M., ROTNITZKY, A. and ZHAO, L. P. (1994). Estimation of regression coefficients when some regressors are not always observed. *J. Amer. Statist. Assoc.* **89** 846–866. MR1294730

SCHUMAKER, L. (1981). *Spline Functions*: *Basic Theory*. Wiley, New York. MR606200

VAN DER VAART, A. W. (1998). *Asymptotic Statistics*. Cambridge Univ. Press. MR1652247

VAN DER VAART, A. W. and WELLNER, J. A. (1996). *Weak Convergence and Empirical Processes*. Springer, New York. MR1385671

ZENG, D. (2004). Estimating marginal survival function by adjusting for dependent censoring using many covariates. *Ann. Statist.* **32** 1533–1555. MR2089133

DEPARTMENT OF BIOSTATISTICS
UNIVERSITY OF NORTH CAROLINA
CHAPEL HILL, NORTH CAROLINA 27599-7420
USA
E-MAIL: dzeng@bios.unc.edu